\documentstyle[12pt]{article}
\begin{document}

\title{Surfaces in Lie sphere geometry and the stationary 
 Davey-Stewartson hierarchy}
\author{{\Large Ferapontov E.V.\thanks{
    Present address:
    Fachbereich Mathematik, SFB 288,
    Technische Universit\"at Berlin,
    10623 Berlin,
    Deutschland,\ \ 
    \hbox{e-mail: {\tt fer@sfb288.math.tu-berlin.de}}}}\\
    Institute for Mathematical Modelling\\
    Academy of Science of Russia, Miusskaya, 4\\
    125047 Moscow, Russia\\
    e-mail: {\tt fer@landau.ac.ru}
}
\date{}
\maketitle

\newtheorem{theorem}{Theorem}

\pagestyle{plain}

\maketitle

\begin{abstract}

We introduce two basic invariant forms which define generic surface in 3-space
uniquely up to Lie sphere equivalence. Two particularly interesting classes
of surfaces associated with these invariants are considered, namely, the
Lie-minimal surfaces and the diagonally-cyclidic surfaces. For diagonally-cyclidic
surfaces we derive the stationary modified Veselov-Novikov equation, whose role in
the theory of these surfaces is similar to that of Calapso's equation in the
theory of isothermic surfaces. Since Calapso's equation itself
turns out to be related to the stationary Davey-Stewartson equation, these
results shed some new light on  differential geometry of the
stationary Davey-Stewartson hierarchy. Diagonally-cyclidic surfaces are the natural
Lie sphere analogs of the isothermally-asymptotic surfaces in projective differential 
geometry for which we also derive the stationary modified Veselov-Novikov equation with 
the different real reduction.

Parallels between invariants of surfaces in Lie sphere geometry and reciprocal
invariants of hydrodynamic  type systems are drawn in the conclusion.

\end{abstract}

\newpage
\section{Introduction}

Lie sphere geometry dates back to the dissertation of Lie in 1872
\cite{Lie}.
After that the subject was extensively developed 
by Blaschke and his coworkers and resulted in publication in 1929 of 
Blaschke's Vorlesungen \"uber Differentialgeometrie \cite{Bla3}, entirely
devoted to the Lie sphere geometry of curves and surfaces.
The modern multidimensional period of the theory was initiated by 
Pinkall's classification of
Dupin hypersurfaces in $E^4$ \cite{Pinkall}, \cite{Pinkall2}. 
We refer also to  Cecil's book \cite{Cecil2}
with the review of the last results in this direction.
Since most of the recent research in Lie sphere geometry is concentrated around
Dupin hypersurfaces and Dupin submanifolds, the general theory of Lie-geometric 
hypersurfaces seems not to be constructed so far.
The aim of this paper is to shed some new light on Lie sphere geometry of (hyper)surfaces
and to reveal its remarkable  applications in the modern theory of integrable systems.

Let $M^2\in E^3$ be a surface in the 3-dimensional space $E^3$ parametrized by the
coordinates $R^1, R^2$ of the lines of curvature. Let $k^1, k^2$ and
$g_{11}{dR^1}^2 + g_{22}{dR^2}^2$ be the principal
curvatures and the metric of $M^2$,
respectively. In sect. 2 we introduce the basic Lie sphere invariants of the surface
$M^2$, namely, the symmetric 2-form
$$
\frac{\partial_1k^1\partial_2k^2}{(k^1-k^2)^2}~dR^1dR^2
\eqno(1.1)
$$
and the conformal class of the cubic form
$$
\partial_1k^1g_{11}{dR^1}^3+\partial_2k^2g_{22}{dR^2}^3,
\eqno(1.2)
$$
$\partial_i=\frac{\partial }{\partial R^i}$, which define "generic"
surface $M^2$ uniquely up to Lie sphere equivalence.
We recall that the group of Lie sphere transformations in $E^{n+1}$ is a contact group, 
generated by conformal transformations and normal shifts, translating each point of the 
surface to a fixed distance $a=const$ along the normal direction.
Conformal transformations and normal shifts generate in $E^{n+1}$ a finite-dimensional
Lie group isomorphic to $SO(n+2,2)$. Lie sphere
transformations can be equivalently characterized as the contact transformations,
mapping spheres into spheres and preserving their oriented contact.
In the implicit form objects (1.1) and (1.2) were contained
already in \cite{Bla3}. Quadratic form (1.1) gives rise to the Lie-invariant
functional
$$
\int \int \frac{\partial_1k^1\partial_2k^2}{(k^1-k^2)^2}~dR^1dR^2
\eqno(1.3)
$$
the extremals of which are known as minimal surfaces in Lie sphere geometry
(\cite{Bla3}, $\S$ 94). Other possible representations of functional (1.3) and its
relationship to the integrable hierarchy of Davey-Stewartson type are discussed in 
sect. 3.  In sect. 4 we investigate the so-called diagonally-cyclidic surfaces
("diagonalzyklidische fl\"achen" in the terminology of \cite{Bla3}, p. 406),  
which can be characterized as surfaces $M^2$  possessing parametrization $R^1, R^2$
by the coordinates of lines of curvature such that the cubic form (1.2)
becomes proportional
to  ${dR^1}^3+{dR^2}^3$. This class of surfaces is a straightforward generalization 
of isothermic surfaces in conformal geometry.
It is demonstrated that for diagonally-cyclidic surfaces the Lie sphere density $U$ defined 
by
$$
U^2=\frac{\partial_1k^1\partial_2k^2}{(k^1-k^2)^2}
$$
satisfies the stationary modified Veselov-Novikov (mVN) equation
$$
\begin{array}{c}
\partial_1^3U+3V\partial_1U+\frac{3}{2}U\partial_1V= 
\partial_2^3U+3W\partial_2U+\frac{3}{2}U\partial_2W
\\
\partial_1W=\partial_2(U^2) \\
\partial_2V=\partial_1(U^2) \\
\end{array}
$$
which in the theory of diagonally-cyclidic surfaces plays a role similar to that of the 
Calapso equation \cite{Calapso}
$$
\partial_1^2\left(\frac{\partial_1\partial_2u}{u}\right)+
\partial_2^2\left(\frac{\partial_1\partial_2u}{u}\right)+
\frac{1}{2}\partial_1\partial_2(u^2)=0
$$
in the theory of isothermic surfaces. Since Calapso's equation itself
turns out to be related to
the stationary Davey-Stewartson (DS) equation (see sect.4), this provides remarkable
differential-geometric interpretation of the stationary flows of  DS hierarchy
and gives new insight into the recent results of \cite{Kon1} --
\cite{Grin} relating DS hierarchy with conformal geometry.
The details of derivation of mVN equation for diagonally-cyclidic surfaces are included 
in the Appendix.

It is quite remarkable that in projective differential geometry there also exists
a class of surfaces (the so-called isothermally-asymptotic surfaces, or $\Phi$-surfaces)
governed by the stationary mVN equation
$$
\begin{array}{c}
p_{xxx}-3Vp_x-\frac{3}{2}pV_x= 
p_{yyy}-3Wp_y-\frac{3}{2}pW_y
\\
W_x=(p^2)_y \\
V_y=(p^2)_x \\
\end{array}
$$
which can be reduced to that presented above by a complex change of variables
$p\to iU, ~~ W\to -W, ~~ V\to -V$. These are surfaces, for which the Darboux cubic form
is proportional to the sum of pure cubes $dx^3+dy^3$  in the appropriate  asymptotic 
parametrization
$x, y$ (particular examples are affine spheres related to  the Tzitzeica equation).
This observation reflects the duality between projective and Lie sphere geometries
due to the Lie's famous line-sphere correspondence (see sect.4).

In sect. 5 we introduce Lie sphere invariants of higher dimensional hypersurfaces
$M^n\in E^{n+1}$, namely the symmetric 2-form
$$
\sum_{i\ne j}\frac{k^i_ik^j_j}{(k^i-k^j)^2}~\omega^i\omega^j
\eqno(1.4)
$$
and the conformal class of the cubic form
$$
\sum_i k^i_ig_{ii}~{\omega^i}^3
\eqno(1.5)
$$
defining "generic" hypersurface uniquely up to Lie sphere equivalence. Here $k^i$ are
principal curvatures, $\omega^i$ are principal covectors,
$\sum g_{ii}{\omega^i}^2$ is the first fundamental form and  coefficients $k_i^i$
are defined by the expansions $dk^i=k^i_j\omega^j$ (we emphasize that
hypersurface $M^n$ of dimension $n\geq 3$ does not necessarily
possess parametrization by the coordinates of lines of curvature). Objects
(1.4) and (1.5) are the Lie-geometric analogs of the second fundamental form
and the Darboux cubic form in projective differential geometry of
hypersurfaces.

In sect. 6-9  interrelations between Lie sphere invariants
and reciprocal invariants of hydrodynamic type systems
$$
u_t^i=v^i_j(u)u^j_x, ~~~ i, j=1,...,n
\eqno(1.6)
$$
are discussed.
We recall that reciprocal transformations are transformations from $x, t$ to the 
new independent variables $X, T$ defined by the formulae
$$
\begin{array}{c}
dX=B(u)dx+A(u)dt\\
dT=N(u)dx+M(u)dt\\
\end{array} 
$$
where $Bdx+Adt$ and $Ndx+Mdt$ are two integrals of system (1.6).
Reciprocal transformations originate from gas dynamics and were extensively investigated
in \cite{Rogers1}, \cite{Rogers2}. In \cite{Fer1}, \cite{Fer2} we introduced reciprocal 
invariants, defining a hydrodynamic type system uniquely up to reciprocal equivalence.
The summary of these results in the 2-component case is given in sect. 6. In sect. 7-8
we recall the necessary information about Hamiltonian systems of hydrodynamic type
and describe the general construction of \cite{Fer5}, \cite{Fer3},
relating Hamiltonian systems (1.6) and hypersurfaces in $E^{n+1}$. The main
property of this correspondence is its "equivariance" in the sense that
Lie sphere transformations of hypersurfaces correspond to "canonical" reciprocal 
transformations, that is, to those reciprocal transformations which preserve
the Hamiltonian structure. In this approach Lie sphere invariants of hypersurfaces
correspond to reciprocal invariants of hydrodynamic type systems, providing thus
their differential-geometric interpretation.

In sect. 9 we write down reciprocal invariants of n-component systems for arbitrary
$n\geq 3$ since they differ from those in case $n=2$.

\section{Invariants of surfaces in Lie sphere geometry}

In \cite{Bla3}, p.392 Blaschke introduced the Lie-invariant 1-forms ~
$\omega^1, ~\omega^2$  ~ ($d\psi, ~d\bar \psi$ in Blaschke's notation)
which assume the following form in the coordinates $R^1, R^2$ of the lines of curvature:
$$
\begin{array}{c}
\omega^1=\frac{k^1_1}{k^1-k^2}\left(
\frac{{k^2_2}^2g_{11}}{{k^1_1}^2g_{22}}\right)^{\frac{1}{6}}dR^1, \\
\ \\
\omega^2=\frac{k^2_2}{k^2-k^1}\left(
\frac{{k^1_1}^2g_{22}}{{k^2_2}^2g_{11}}\right)^{\frac{1}{6}}dR^2. \\
\end{array}
\eqno(2.1)
$$

{\bf Remark 1.} In order to check Lie-sphere invariance of 
the 1-forms $\omega^1, \omega^2$ it is sufficient to check their invariance under
the inversions and normal shifts, which can be verified by a direct calculation.
Moreover, forms (2.1) do not change if the principal curvatures $k^i$ are
replaced by the radii of principal curvatures $w^i=\frac{1}{k^i}$ and the first
fundamental form $g_{ii}$ by the third fundamental form $G_{ii}={k^i}^2g_{ii}$.

{\bf Remark 2.} Similar invariant 1-forms arise in the M\"obius (conformal) and
Laguerre geometries, which are subcases of the Lie sphere geometry.
In the M\"obius geometry we have the invariant 1-forms
$$
\frac{\partial_1k^1}{k^1-k^2}dR^1, ~~~~~ 
\frac{\partial_2k^2}{k^2-k^1}dR^2 
$$ 
and the invariant quadratic form
$$
(k^1-k^2)^2(g_{11}{dR^1}^2+g_{22}{dR^2}^2),
$$
while in the Laguerre geometry they are
$$
\frac{\partial_1w^1}{w^1-w^2}dR^1, ~~~~~ \frac{\partial_2w^2}{w^2-w^1}dR^2 
$$ 
and 
$$
(w^1-w^2)^2(G_{11}{dR^1}^2+G_{22}{dR^2}^2),
$$
respectively (see \cite{Bla3}, \cite{Wang}, \cite{AkGold1}).

As far as $\omega^1$ and $\omega^2$ are invariant under Lie sphere 
transformations, so do the quadratic form
$$
-\omega^1\omega^2=\frac{\partial_1k^1\partial_2k^2}{(k^1-k^2)^2}~dR^1dR^2
$$
and the qubic form
$$
{\omega^1}^3-{\omega^2}^3=
\frac{\partial_1k^1\partial_2k^2}{(k^1-k^2)^3\sqrt{g_{11}g_{22}}}
(\partial_1k^1g_{11}{dR^1}^3+\partial_2k^2g_{22}{dR^2}^3),
\eqno(2.2)
$$
giving rise to (1.1) and (1.2), respectively. 
The reason for introducing these objects is their additional
invariance under the interchange of indices 1 and 2 so that they play a role
similar to that of "symmetric functions of the roots of polynomial" in the Viete 
theorem. Hence they define tensors which can be effectively computed 
in an arbitrary coordinate system (without solving algebraic equations).

In \cite{Bla3}, $\S 85$ it is proved that up to certain exeptional cases
a generic surface in $E^3$ is determined by the corresponding 1-forms 
$\omega^1, \omega^2$ uniquely up to Lie sphere transformations. Since 
we can reconstruct $\omega^1, \omega^2$ from the quadratic form (1.1) and the 
conformal class of the cubic form (1.2) (the multiple in (2.2) is not
essential), we can formulate the following

{\bf Theorem 1.}
{\it A generic surface $M^2\in R^3$ is defined by the quadratic form 
$$
\frac{\partial_1k^1\partial_2k^2}{(k^1-k^2)^2}~dR^1dR^2
$$
and the conformal class of the cubic form
$$
\partial_1k^1g_{11}{dR^1}^3+\partial_2k^2g_{22}{dR^2}^3
$$
uniquely up to Lie sphere transformations}.

The vanishing of the qubic form is equivalent to the conditions
$\partial_1k^1=\partial_2k^2=0$
which specify the so-called cyclids of Dupin. We recall that the vanishing of 
Darboux's cubic form in projective geometry  specifies quadrics, which 
are thus projective "duals" of cyclids of Dupin.

{\bf Remark}. The principal directions of the surface $M^2$ can be characterized as
the zero directions of  quadratic form (1.1). On the other hand, they are exactly  
those directions, where cubic form (1.2) reduces to the sum of pure cubes 
(without mixed terms). It should be  pointed out, that
any cubic form on the plane can be reduced to the sum of cubes, and the directions
where it assumes the desired form are defined uniquely.
In projective differential geometry of surfaces these are asymptotic directions of
the Darboux cubic form.

\section{Minimal surfaces in Lie sphere geometry}

Lie-minimal surfaces are defined as the extremals of Lie-invariant functional (1.3)
$$
\int \int \frac{\partial_1k^1\partial_2k^2}{(k^1-k^2)^2}~dR^1dR^2
$$
which  is a natural analog of the Willmore  functional
$$
\int \int (k^1-k^2)^2\sqrt {g_{11}g_{22}}~dR^1dR^2
$$
in conformal geometry and the invariant functional
$$
\int \int (w^1-w^2)^2\sqrt {G_{11}G_{22}}~dR^1dR^2
$$
in the Laguerre geometry. Due to the obvious identity
$$
\frac{\partial_1k^1\partial_2k^2}{(k^1-k^2)^2}~dR^1\wedge dR^2=
\frac{\partial_2k^1\partial_1k^2}{(k^1-k^2)^2}~dR^1\wedge dR^2-
d\left(\frac{dk^2}{k^1-k^2}\right)
$$
we see that for compact surfaces with $k^1\ne k^2$ (for instance, immersed tori)
functional (1.3) coincides with the functional
$$
\int \int \frac{\partial_2k^1\partial_1k^2}{(k^1-k^2)^2}~dR^1dR^2=
-\int \int ab~ dR^1dR^2,
\eqno(3.1)
$$
where we introduced the notation
$a=\frac{\partial_2k^1}{k^2-k^1}, ~~ b=\frac{\partial_1k^2}{k^1-k^2}$. We recall that in 
terms of the coefficients $a$ and $b$ the  Peterson-Codazzi equations 
of the surface $M^2$ can be written as follows
$$
\partial_2\ln \sqrt {g_{11}}=a, ~~~~ \partial_1\ln \sqrt {g_{22}}=b
$$
while the equation for  the radius-vector $\vec r$ assumes the form 
$$
\partial_1\partial_2\vec r=a\partial_1\vec r+b\partial_2\vec r
\eqno(3.2)
$$
manifesting the fact that the net of the lines of curvature is conjugate.

Introducing the rotation coefficients $\beta_{12}, \beta_{21}$ by the formulae
$$
\beta_{12}=\frac{\partial_1\sqrt{g_{22}}}{\sqrt{g_{11}}}=
b~\frac{\sqrt{g_{22}}}{\sqrt{g_{11}}}, ~~~~~~
\beta_{21}=\frac{\partial_2\sqrt{g_{11}}}{\sqrt{g_{22}}}=
a~\frac{\sqrt{g_{11}}}{\sqrt{g_{22}}},
$$
we can rewrite our functional in the form
$$
\int \int  \beta_{12} \beta_{21} ~ dR^1dR^2
$$
which has the meaning of the "integral squared rotation" of the surface $M^2$.

Written in the form (3.1) our functional is closely related to the simplest 
quadratic conservation law of the (2+1)-dimensional hierarchy of the
Davey-Stewartson (DS) type. To make it clear we recall the construction of \cite{Kon3},
which defines the DS-type flow on conjugate nets in $E^3$. Let $M^2$ be a surface
parametrized by conjugate coordinates $R^1, R^2$. The radius-vector $\vec r$ of                              
such surface satisfies the equation (3.2) 
(at the moment we do not assume that our conjugate
net is the net of lines of curvature). Let us define  evolution of $M^2$
with respect to the "time" $t$ by the formula
$$
\vec r_t=\alpha{\partial_1}^2\vec r + \beta {\partial_2}^2\vec r
         +p\partial_1\vec r+q\partial_2\vec r,
\eqno(3.3)
$$
where $\alpha, \beta =const$. The compatibility conditions of (3.2) and (3.3) give
rise to the integrable system
$$
\begin{array}{c}
a_t=\beta{\partial_2}^2a-\alpha {\partial_1}^2a+
2\alpha \partial_1(ab)+\beta \partial_2(a^2)+p\partial_1a+\partial_2(qa),\\
\ \\
b_t=\alpha{\partial_1}^2b-\beta {\partial_2}^2b+
2\beta \partial_2(ab)+\alpha \partial_1(b^2)+q\partial_2b+\partial_1(pb),\\
\ \\
\partial_2p+2\alpha \partial_1a=0,\\
\ \\
\partial_1q+2\beta \partial_2b=0,\\
\end{array}
\eqno(3.4)
$$
whose relationship to the DS system was clarified in \cite{Kon3}. In fact system (3.4)
is a linear combination of two simpler systems, corresponding to the choices
$(\alpha =1, \beta =0, q=0)$ and $(\alpha =0, \beta =1, p=0)$, namely
$$
\begin{array}{c}
a_t=-{\partial_1}^2a+
2\partial_1(ab)+p\partial_1a,\\
\ \\
b_t={\partial_1}^2b+
\partial_1(b^2)+\partial_1(pb),\\
\ \\
\partial_2p+2\partial_1a=0,\\
\end{array}
\eqno(3.5)
$$
and
$$
\begin{array}{c}
a_t={\partial_2}^2a
+ \partial_2(a^2)+\partial_2(qa),\\
\ \\
b_t=- {\partial_2}^2b+
2\partial_2(ab)+q\partial_2b,\\
\ \\
\partial_1q+2\partial_2b=0,\\
\end{array}
\eqno(3.6)
$$
respectively. Both systems (3.5) and (3.6) commute according to 
the general discussion in \cite{Mikh} and possess the quadratic integral
$$
\int \int ab~ dR^1dR^2
$$
coinciding with (3.1). We emphasize, however, that the nets of lines of curvature
are not preserved (in general) by these $t$-evolutions,
although always remain conjugate.

\section{Cyclidic curves and diagonally-cyclidic surfaces}

With any surface $M^2$ we associate a 3-web of curves
(that is, three 1-parameter families of curves) formed by the lines of
curvature and cyclidic curves ("zyklidische kurven" in the terminology of Blaschke
\cite{Bla3}, \S 86)
which are the zero directions of cubic form (1.2). In view
of formula (2.2) the curves of this 3-web can be defined in terms of  1-forms (2.1)
by the equations
$$
\omega^1=0, ~~~~~ \omega^2=0, ~~~~~ \omega^1-\omega^2=0,
\eqno(4.1)
$$
respectively. Cyclidic curves naturally arise in the attempt to find those cyclids of
Dupin, which are the "best" tangents to a given surface $M^2$ at a given point
(see \cite{Bla3}, \S 86 for the details).
These curves are the natural analogs of the Darboux curves in projective differential
geometry. Let us compute the connection 1-form $\omega$ of the 3-web (4.1), that is,
the 1-form which is uniquely determined by the equations
$$
d\omega^1=\omega \wedge \omega^1, ~~~
d\omega^2=\omega \wedge \omega^2,
$$
(see \cite{Bla1}, \cite{Bla2} for the introduction in web geometry).
A direct computation results in
$$
\begin{array}{c}
\omega =
       \frac{1}{3}\left(\frac{\partial_1\partial_2k^2}{\partial_2k^2}+
             \frac{\partial_1k^1}{k^1-k^2}\right)dR^1+
       \frac{1}{3}\left(\frac{\partial_1\partial_2k^1}{\partial_1k^1}+
             \frac{\partial_2k^2}{k^2-k^1}\right)dR^2 +\\
\ \\
\frac{1}{3}d\ln \frac{\partial_1k^1\partial_2k^2}
{(k^1-k^2)^5\sqrt{g_{11}g_{22}}}.
\end{array}
\eqno(4.2)             
$$
Since both $\omega^1, \omega^2$ are invariant under Lie sphere transformations, 
so does the connection 1-form $\omega$. From (4.2) it immediately follows that
the curvature form $d\omega$ of 3-web (4.1) is given by
$$
d\omega=\frac{1}{3}d\Omega,
$$
where
$$
\Omega =
       \left(\frac{\partial_1\partial_2k^2}{\partial_2k^2}+
             \frac{\partial_1k^1}{k^1-k^2}\right)dR^1+
       \left(\frac{\partial_1\partial_2k^1}{\partial_1k^1}+
             \frac{\partial_2k^2}{k^2-k^1}\right)dR^2.
\eqno(4.3)             
$$
As we will see in sect. 6 the object analogous to (4.3) arises in the theory of
reciprocal invariants of hydrodynamic type systems.

An interesting class of diagonally-cyclidic surfaces ("diagonalzyklidische fl\"achen"
in the terminology of \cite{Bla3}, p.406) is specified by the requirement, that
3-web (4.1) is hexagonal or, equivalently, has zero curvature:
$$
d\omega=\frac{1}{3}d\Omega=0.
$$
In this case there exist coordinates $R^1, R^2$ along the lines of curvature
(note, that we have a reparametrization
freedom $R^i\to \varphi^i(R^i)$), where $\omega^1, \omega^2$ assume the form
$$
\omega^1=pdR^1, ~~~ 
\omega^2=-pdR^2
$$
with nonzero common multiple $p$. Since in these coordinates the cubic form
${\omega^1}^3-{\omega^2}^3$ is proportional to ${dR^1}^3+{dR^2}^3$, these surfaces 
are the natural Lie-sphere analogs of isothermic surfaces in 
conformal geometry (we emphasize, that the class of isothermic surfaces is not 
invariant under the full group of Lie sphere transformations).

Equations governing diagonally-cyclidic surfaces can be easily written down as follows.
Since the cubic form
$$
\partial_1k^1 g_{11}{dR^1}^3+\partial_2k^2 g_{22}{dR^2}^3
$$
is proportional to
$$
{dR^1}^3+{dR^2}^3,
$$
we can put
$$
g_{11}=\frac{e^{2\rho}}{\partial_1k^1}, ~~~~
g_{22}=\frac{e^{2\rho}}{\partial_2k^2}.
$$
Inserting this representation in the Gauss-Peterson-Codazzi equations,
we arrive at the following system for $k^1, k^2, \rho$:
$$
\begin{array}{c}
\partial_1\rho = b+\frac{1}{2}\frac{\partial_1\partial_2k^2}{\partial_2k^2}, \\
\ \\
\partial_2\rho = a+\frac{1}{2}\frac{\partial_1\partial_2k^1}{\partial_1k^1}, \\
\ \\
\partial_1k^1\left(\partial_1b+\frac{b}{2}\partial_1\ln
{\frac{\partial_1k^1}{\partial_2k^2}}\right)+
\partial_2k^2\left(\partial_2a+\frac{a}{2}\partial_2\ln
{\frac{\partial_2k^2}{\partial_1k^1}}\right)+ 
k^1k^2e^{2\rho}=0,
\end{array}
\eqno(4.4)
$$
where
$a=\frac{\partial_2k^1}{k^2-k^1}, ~~
 b=\frac{\partial_1k^2}{k^1-k^2}.$

{\bf Remark}. System (4.4) is a Lie-sphere analog of the system
$$
\begin{array}{c}
\partial_1\rho=b,\\
\partial_2\rho=a,\\
\partial_1^2\rho+\partial_2^2\rho+k^1k^2e^{2\rho}=0,
\end{array}
\eqno(4.5)
$$
describing isothermic surfaces in conformal geometry. The integrability and
 discretization (see \cite{Sym}, \cite{Bobenko}) of system
(4.5) are based on the $SO(4,1)$-linear problem, which comes from the following
geometric fact: all isothermic surfaces possess Ribacour transformations
preserving the metric up to a conformal factor. Moreover, the spectral parameter
is due to the following scailing symmetry of system (4.5):
$$
\tilde R^1=\frac{1}{c}R^1, ~~~ \tilde R^2=\frac{1}{c}R^2,
~~~ \tilde k^1=ck^1, ~~~ \tilde k^2=ck^2, ~~~ c=const.
$$
We recall also that system (4.5) can be 
rewritten as the single fourth-order Calapso equation \cite{Calapso}
$$
\partial_1^2\left(\frac{\partial_1\partial_2 u}{u}\right)+
\partial_2^2\left(\frac{\partial_1\partial_2 u}{u}\right)+
\frac{1}{2}\partial_1\partial_2 (u^2)=0
\eqno(4.6)
$$
for the conformal factor  $u=e^{\rho}(k^1-k^2)$.
Another possible approach to the integrability of isothermic surfaces is based on the
relationship of the Calapso equation to the DS-II equation
$$
\begin{array}{c}
iu_t+u_{xx}-u_{yy}+uv=0,\\
v_{xx}+v_{yy}=\vert u\vert ^2_{xx}- \vert u\vert ^2_{yy}
\end{array}
$$
which in the stationary case assumes the form
$$
\begin{array}{c}
u_{xx}-u_{yy}+uv=0,\\
v_{xx}+v_{yy}=\vert u\vert ^2_{xx}- \vert u\vert ^2_{yy}.
\end{array}
$$
Excluding $v$ we arrive at the fourth-order equation with respect to $u$
$$
\triangle \left(\frac{u_{xx}-u_{yy}}{u}\right)+
\vert u\vert ^2_{xx}- \vert u\vert ^2_{yy}=0
$$
coinciding with (4.6) after the transformation $R^1=x+y, ~ R^2=x-y$ and
the reduction $u=\bar u$.

\bigskip

Similar approaches can be applied to system (4.4). 
Indeed, all diagonally-cyclidic surfaces possess
Ribacour transformations preserving the cubic form up to a conformal factor
(\cite{Bla3}, p.420).
Moreover, system (4.4) possesses a similar scaling symmetry
$$
\tilde R^1=\frac{1}{c}R^1, ~~~ \tilde R^2=\frac{1}{c}R^2,
~~~ \tilde k^1=c^3k^1, ~~~ \tilde k^2=c^3k^2,
$$
which is responsible for the spectral parameter. We hope to develop this geometric 
approach elsewhere.  Another approach to the integrability of diagonally-cyclidic
surfaces is based on their remarkable relationship to the mVN equation
$$
\begin{array}{c}
U_t=U_{xxx}-U_{yyy}+3U_xV-3U_yW +\frac{3}{2}UV_x-\frac{3}{2}UW_y\\

W_x=(U^2)_y \\
V_y=(U^2)_x \\
\end{array}
$$
introduced in \cite{Bogdanov} which is the third-order flow in the
DS-I hierarchy.
In the Appendix we demonstrate, that the Lie sphere
density $U$ defined by
$$
U^2=\frac{\partial_1k^1\partial_2k^2}{(k^1-k^2)^2}
$$
satisfies in case of diagonally-cyclidic surfaces the stationary mVN equation
$$
\begin{array}{c}
U_{xxx}+3U_xV+\frac{3}{2}UV_x=U_{yyy}+3U_yW+\frac{3}{2}UW_y \\

W_x=(U^2)_y \\
V_y=(U^2)_x \\
\end{array}
$$
(here $\partial_1=\partial_x, ~ \partial_2=\partial_y$). Note that $U^2$ is a conserved
density of the mVN equation. A passage from equations (4.4) to he mVN equation
requires quite complicated calculations which were performed with Mathematica
(see the Appendix).

Particular solutions of the stationary mVN equation can be obtained by the ansatz
$$
W=-\frac{2}{3}\frac{U_{yy}}{U}+\frac{1}{3}\left(\frac{U_{y}}{U}\right)^2
$$
$$
V=-\frac{2}{3}\frac{U_{xx}}{U}+\frac{1}{3}\left(\frac{U_{x}}{U}\right)^2
$$
where $U$ satisfies the Tzitzeica equation
$$
(\ln U)_{xy}=-U^2+\frac{c}{U}, ~~~ c=const.
$$
In this case the first equation is satisfied identically, while the last two 
become just two conservation laws of the Tzitzeica equation.

We have demonstrated that the stationary flows of the DS hierarchy
have a natural interpretation within the contexts of 
conformal and Lie sphere geometries. This observation agrees with the results of
\cite{Kon1}, \cite{Kon2}, \cite{Taiman3}, \cite{Taiman1}, \cite{Taiman2}, \cite{Grin}
where dynamics of surfaces, induced by the odd flows of
DS-II hierarchy, was investigated. In particular, it was argued that
the integrals of DS-II hierarchy define conformally invariant functionals
and the corresponding stationary flows define certain conformally invariant classes
of surfaces. We hope that results presented above contribute to these
investigations. We emphasize also, that stationary points of some of the
DS-integrals are probably invariant under the full group of Lie sphere transformations,
rather then just  conformal group (for example, the class of diagonally-cyclidic 
surfaces, corresponding to the stationary mVN equation, is invariant under the
full Lie sphere group).

{\bf Remark.} Diagonally-cyclidic
surfaces have a natural projective counterpart, namely,
the so-called isothermally-asymptotic surfaces for which the 3-web, 
formed by the asymptotic
lines and Darboux's curves is hexagonal
(Darboux's curves  are the zero curves of Darboux's cubic form).
This class of surfaces can be equivalently characterized by the existence 
of asymptotic coordinates where Darboux's cubic form becomes isothermic.
We refer to \cite{Finikov} for further discussion and references conserning
isothermally-asymptotic surfaces ($\Phi$-surfaces in the terminology of \cite{Finikov}).
Isothermally-asymptotic surfaces are related to a different real reduction of the 
mVN equation. Here we present the details of it's derivation. Let $M^2$ be a surface 
in projective space parametrized by asymptotic coordinates $x, y$ with the radius-vector 
$\vec r$ satisfying the equations
$$
\begin{array}{c}
\vec r_{xx}=a\vec r_x+p\vec r_y, \\
\vec r_{yy}=q\vec r_x+b\vec r_y.
\end{array}
\eqno(4.7)
$$
With any surface (4.7) we associate
the symmetric 2-form
$$
pq~dxdy
\eqno(4.8)
$$
and the conformal class of the Darboux cubic form
$$
pdx^3+qdy^3
\eqno(4.9)
$$
which define "generic" surface $M^2$ uniquely up to projective equivalence and play a 
role similar to that of (1.1) and (1.2) in the Lie sphere geometry. Darboux's curves 
are the zero curves of the Darboux cubic form. The hexagonality conditions of the 3-web
formed by the asymptotic lines and the Darboux curves is equivalent to the existence 
of asymptotic parametrization $x, y$ such that cubic form (4.9) becomes proportional to 
$$
dx^3+dy^3,
$$
that is, to the condition $p=q$. An important subclass of isothermally-asymptotic 
surfaces are the proper affine spheres, for which the radius-vector $\vec r$ satisfies 
the equations
$$
\begin{array}{c}
\vec r_{xx}=-\frac{p_x}{p}\vec r_x+p\vec r_y \\
\ \\
\vec r_{yy}=p\vec r_x-\frac{p_y}{p}\vec r_y \\
\ \\
\vec r_{xy}=\frac{1}{p}\vec r
\end{array}
$$
with $p$ satisfying the Tzitzeica equation
$$
(\ln p)_{xy}=p^2-\frac{1}{p}.
$$
In the case $p=q$ the compatibility conditions of (4.7) reduce to
$$
\begin{array}{c}
(p_x+ap+\frac{1}{2}b^2-b_y)_x=\frac{3}{2}(p^2)_y \\
\ \\
(p_y+bp+\frac{1}{2}a^2-a_x)_y=\frac{3}{2}(p^2)_x \\
\ \\
a_y=b_x.
\end{array}
\eqno(4.10)
$$
Equations (4.10) can be rewritten in the form
$$
\begin{array}{c}
b_y=p_x+ap+\frac{1}{2}b^2-\frac{3}{2}W, ~~~~~~ W_x=(p^2)_y \\
\ \\
a_x=p_y+bp+\frac{1}{2}a^2-\frac{3}{2}V, ~~~~~~ V_y=(p^2)_x \\
\ \\
a_y=f, ~~~~~ b_x=f.
\end{array}
\eqno(4.11)
$$
Crossdifferentiation of (4.11) gives the expressions for $f_x, f_y$
$$
\begin{array}{c}
f_x=p_{yy}+ap^2+\frac{1}{2}pb^2-\frac{3}{2}pW+bp_y+af-(p^2)_x \\
\ \\
f_y=p_{xx}+bp^2+\frac{1}{2}pa^2-\frac{3}{2}pV+ap_x+bf-(p^2)_y \\
\end{array}
\eqno(4.12)
$$
the compatibility conditions of which result in the stationary mVN equation
$$
\begin{array}{c}
p_{xxx}-3Vp_x-\frac{3}{2}pV_x= 
p_{yyy}-3Wp_y-\frac{3}{2}pW_y
\\
W_x=(p^2)_y \\
V_y=(p^2)_x \\
\end{array}
\eqno(4.13)
$$
which can be reduced to that presented above by a complex change of variables
$p\to iU, ~~ W\to -W, ~~ V\to -V$. Integrating the compatible systems (4.12),
(4.11) and (4.7) for a given solution $p, W, V$ of (4.13) we arrive at the 
explicit formula for the radius-vector $\vec r$. 

Particular solutions of the stationary mVN equation  (4.13) 
can be obtained by the ansatz
$$
W=\frac{2}{3}\frac{p_{yy}}{p}-\frac{1}{3}\left(\frac{p_{y}}{p}\right)^2
$$
$$
V=\frac{2}{3}\frac{p_{xx}}{p}-\frac{1}{3}\left(\frac{p_{x}}{p}\right)^2
$$
where $p$ satisfies the Tzitzeica equation
$$
(\ln p)_{xy}=p^2+\frac{c}{p}, ~~~ c=const.
$$
In this case the first equation is satisfied identically, while the last two 
become just two conservation laws of the Tzitzeica equation.
If $c\ne 0$, then it can be normalized to $-1$; the corresponding surfaces 
are the proper affine spheres. The case $c=0$ 
corresponds to the improper affine spheres 
whose affine normals are parallel and the radius-vector $\vec r$
satisfies the equations
$$
\begin{array}{c}
\vec r_{xx}=-\frac{p_x}{p}\vec r_x+p\vec r_y \\
\ \\
\vec r_{yy}=p\vec r_x-\frac{p_y}{p}\vec r_y \\
\ \\
\vec r_{xy}=\frac{1}{p}\vec l
\end{array}
$$
where $\vec l$ is a constant vector (direction of the affine normal) and $p$
satisfies the Liouville equation
$$
(\ln p)_{xy}=p^2.
$$
The fact that Tzitzeica's equation defines a subclass of solutions of the 
stationary VN equation was observed recently in \cite{Kon4}
(I would like to thank W.~Schief for providing me with this reference).  
Our results give differential-geometric interpretation of this formal observation. 

Isothermally-asymptotic surfaces are known to possess B\"acklund transformations
such that the initial and the transformed surfaces are two focal surfaces of a 
W-congruence, which preserves asymptotic lines and the Darboux curves. This can 
be the starting point for the modern approach to isothermally-asymptotic surfaces,
their discretization in the spirit of \cite{Bobenko2}, etc.

\section{Invariants of higher dimensional hypersurfaces in Lie sphere geometry}

In this section we announce several results on Lie sphere geometry on higher 
dimensional hypersurfaces, postponing the detailed proofs to a separate publication.

Let $M^n$ be hypersurface with principal curvatures $k^i$ and principal 
covectors
$\omega^i$, so that the $i$-th principal direction of $M^n$ is defined by the
equations $\omega^j=0,~~ j\ne i$. It should be pointed out that generic hypersurface
of dimension $\geq 3$ does not possess parametrization $R^i$ by the lines of 
curvature as in the 2-dimensional case. 
Differentiating covectors $\omega^i$ and  principal curvatures $k^i$ we arrive at the 
structure equations
$$
d\omega^i=c^i_{jk}\omega^j \wedge \omega^k,
\eqno(5.1)
$$
$$
dk^i=k^i_j\omega^j.
$$
Let also 
$$
ds^2=\sum_{1}^{n}g_{ii}{\omega^i}^2
$$
be the first fundamental form of hypersurface $M^n$.

{\bf Theorem 2.} {\it A generic hypersurface $M^n$ $(n\geq 3)$ is defined by 
quadratic form (1.4)
$$
\sum_{i\ne j}\frac{k^i_ik^j_j}{(k^i-k^j)^2}~\omega^i\omega^j
$$
and the conformal class of  cubic form (1.5)
$$
\sum_i k^i_ig_{ii}~{\omega^i}^3
$$
uniquely up to Lie sphere equivalence.}

As "generic" it is sufficient to understand a surface with $k^i_i\ne 0$.
This genericity assumption is essential since, for instance, there exist examples
of Dupin hypersurfaces (that is, hypersurfaces with $k^i_i=0$) which are not 
Lie-equivalent. Theorem 2 is an analog of the corresponding theorem in projective 
differential geometry stating that hypersurface in projective space $P^n$ of dimension
$n\geq 4$ is uniquely determined by the conformal classes of of its second fundamental 
form and Darboux's cubic form  (see \cite{AkGold2} for the exact statements and 
further references). 

{\bf Remark 1}. In case $k^i_i\ne 0$ cubic form (1.5) encodes all the information 
about the lines of curvature of hypersurface $M^n$. Indeed, principal directions 
are uniquely defined as those directions where cubic form (1.5) reduces to the
sum of pure cubes (without mixed terms). 
Moreover, principal directions are zero directions of quadratic form (1.4).
However, this last condition does not define them uniquely as in the 
2-dimensional situation.

{\bf Remark 2}. The invariant quadratic form (1.4) defines the invariant
volume form, giving rise  in the case $n=3$ to the invariant functional
$$
\int \int \int \frac{k^1_1k^2_2k^3_3}{(k^1-k^2)(k^1-k^3)(k^2-k^3)}
~\omega^1\omega^2\omega^3,
\eqno(5.2)
$$
the extremals of which should be called minimal hypersurfaces in Lie sphere geometry
in analogy with the 2-dimensional case.
It does not look likely that this functional was investigated so far.

{\bf Remark 3}. 
In principle for $n\geq 3$ there exist additional Lie-sphere invariants besides
those mentioned in Theorem 2, namely:

1. The cross-ratios
$$
\frac{(k^i-k^j)(k^n-k^l)}
     {(k^n-k^j)(k^i-k^l)}
\eqno(5.3)
$$
of any four principal curvatures.

2. The covectors
$$
\frac{k^i_i(k^j-k^l)}{(k^i-k^j)(k^i-k^l)}
~\omega^i ~~~~
(i\ne j\ne l).
\eqno(5.4)
$$
For instance, in case $n=3$ we have three invariant covectors
$$
\Omega^1=\frac{k^1_1(k^2-k^3)}{(k^1-k^2)(k^1-k^3)}~\omega^1, ~~~
\Omega^2=\frac{k^2_2(k^3-k^1)}{(k^2-k^1)(k^2-k^3)}~\omega^2, 
$$
$$
\Omega^3=\frac{k^3_3(k^1-k^2)}{(k^3-k^1)(k^3-k^2)}~\omega^3
$$
giving rise to the invariant quadratic form
$$
{\Omega^1}^2+{\Omega^2}^2+{\Omega^3}^2
$$
whose volume functional
$$
\int \int \int \Omega^1\Omega^2\Omega^3
$$
coincides with (5.2).

3. Conformal class of the quadratic form
$$
g_{11}\left(\prod_{l\ne 1}(k^1-k^l)\right)^{\frac{2}{n-2}}{\omega^1}^2
+\ldots +
g_{nn}\left(\prod_{l\ne n}(k^n-k^l)\right)^{\frac{2}{n-2}}{\omega^n}^2.
\eqno(5.5)
$$
Lie-invariant class of hypersurfaces with conformally flat quadratic form (5.5)
deserves a special investigation.

4. Differential $d\Omega$ of the 1-form
$$
\Omega=\left(\sum_{l\ne 1}\frac{k^l_1-\frac{k^1_1}{n-1}}{k^1-k^l}
         \right)\omega^1+\ldots +
       \left(\sum_{l\ne n}\frac{k^l_n-\frac{k^n_n}{n-1}}{k^n-k^l}
         \right)\omega^n.
\eqno(5.6)
$$

In generic case $k^i_i\ne 0$ objects (5.3)-(5.6) can be expressed 
through the forms (1.4) and (1.5). However they are important in the nongeneric 
situations, when some (or all) of $k^i_i$ vanish so that  
(1.4) and (1.5) become zero. In particular, cross-ratios of
principal curvatures play essential role in the study of Dupin hypersurfaces -- see
\cite{Cecil1}, \cite{Miyaoka}. In this respect it seems 
interesting to understand the role of
conformal class (5.5) and the 2-form $d\Omega$ in the modern 
Lie-geometric approach to Dupin hypersurfaces.

For hypersurfaces with nonholonomic nets of lines of curvature Theorem 2 leads to a nice 
geometric corollary, which we will discuss in the simplest nontrivial 3-dimensional case.
Let us consider the structure equations (5.1) of the 3-dimensional hypersurface $M^3$.
Then there are only two possibilities:

1. Holonomic case: all three coefficients $c^1_{23}, c^2_{31}, c^3_{12}$ are equal 
to zero. Such hypersurfaces possess parametrization by the lines of curvature.

2. Nonholonomic case: all three coefficients
$c^1_{23}, c^2_{31}, c^3_{12}$ are nonzero. 

It immediately follows from the Peterson-Codazzi equations that for $n=3$
intermediate cases are forbidden. In the nonholonomic case we can normalize
covectors $\omega^1, \omega^2, \omega^3$ in such a way that the structure equations
assume the form
$$
\begin{array}{c}
d\omega^1=a\omega^1\wedge \omega^2+b\omega^1\wedge \omega^3+\omega^2\wedge \omega^3,\\
d\omega^2=p\omega^2\wedge \omega^1+q\omega^2\wedge \omega^3+\omega^3\wedge \omega^1,\\
d\omega^3=r\omega^3\wedge \omega^1+s\omega^3\wedge \omega^2+\omega^1\wedge \omega^2.\\
\end{array}
\eqno(5.7)
$$
This normalization fixes $\omega^i$ uniquely. As follows from the results of \cite{Fer6},
in the 3-dimensional nonholonomic situation Peterson-Codazzi equations completely 
determine quadratic form (1.4) and qubic form (1.5) through the coefficients
$a, b, p, q, r, s$ in the structure equations (5.7). Hence we can formulate the 
following result.

{\bf Theorem 3}. {\it Nonholonomic 3-dimensional hypersurface $M^3$ is defined by its
structure equations (5.7) uniquely up to Lie sphere equivalence.}

We can reformulate this result as follows: two 3-dimensional nonholonomic hypersurfaces
are Lie-equivalent if and only if there exists a point correspondence between them,
mapping the lines of curvature of one of them onto the lines of curvature of the other.
Hence 3-dimensional nonholonomic hypersurface is
uniquely determined by geometry of it's 
lines of curvature.

This Theorem should remain valid for higher-dimensional hypersurfaces
if we generalize the notion of "nonholonomicity" in a proper way (e.g. $c^i_{jk}\ne 0$
for all $i\ne j\ne k$ which probably can be weakend).

\section{Reciprocal transformations of hydrodynamic type systems.
          Reciprocal invariants}

In this section we consider 2-component systems of hydrodynamic type
$$
u_t^i=v^i_j(u)u^j_x, ~~~~~ i, j=1, 2
\eqno(6.1)
$$
which naturally arise in polytropic gas dynamics, chromatography,
plasticity , etc. and describe wide variety of models of continuous media.
The main advantage of the 2-component case is the existence of the so-called Riemann 
invariants: coordinates, where equations (6.1) assume the diagonal form
$$
\begin{array}{c}
R_t^1=\lambda^1(R)R^1_x,\\
R_t^2=\lambda^2(R)R^2_x,\\
\end{array} 
\eqno(6.2)
$$
considerably simplifying their investigation.
Any system (6.2) possesses infinitely many conservation laws
$$
h(R)dx+g(R)dt
\eqno(6.3)
$$
with the densities $h(R)$ and the fluxes $g(R)$ governed by the equations
$$
\partial_ig=\lambda^i\partial_ih, ~~~~~ i=1, 2
\eqno(6.4)
$$
($\partial_i={\partial}/{\partial R^i}$)
which are completely equivalent to the condition $h_t=g_x$, manifesting closedness 
of  1-form (6.3). Crossdifferentiation of (6.4) results in the
second-order equation
$$
\partial_1\partial_2h=\frac{\partial_2\lambda^1}{\lambda^2-\lambda^1}\partial_1h+
\frac{\partial_1\lambda^2}{\lambda^1-\lambda^2}\partial_2h, 
\eqno(6.5)
$$
for the conserved densities of system (6.2). Thus conservation laws of system (6.2) 
depend on two arbitrary functions of 
one variable. Let us choose two particular conservation laws
$B(R)dx+A(R)dt$ and $N(R)dx+M(R)dt$ and introduce new independent variables
$X, T$ by the formulae
$$
\begin{array}{c}
dX=Bdx+Adt\\
dT=Ndx+Mdt\\
\end{array} 
\eqno(6.6)
$$
which are correct since the right hand sides are closed. Changing from $x, t$ to 
$X, T$ in (6.2) we arrive at the transformed system
$$
\begin{array}{c}
R_t^1=\Lambda^1(R)R^1_x\\
R_t^2=\Lambda^2(R)R^2_x\\
\end{array} 
\eqno(6.7)
$$
where the new characteristic velocities $\Lambda^i$ are given by the formulae
$$
\Lambda^i=\frac{\lambda^iB-A}{M-\lambda^iN}, ~~~~~ i=1, 2.
\eqno(6.8)
$$
{\bf Remark}. In principle one can apply transformation (6.6) directly
to system (6.1) without rewriting it in Riemann invariants. In this case the transformed 
equations assume the form
$$
u_T^i=V^i_j(u)u^j_X,
$$
with the new matrix $V$ given by
$$
V=(Bv-AE)(ME-Nv)^{-1}, ~~~ E=id.
$$

Transformations of  type (6.6) are known as "reciprocal" and have been 
extensively investigated in \cite{Rogers1}, \cite{Rogers2}
(see also \cite{Serre1} and \cite{Fer1} - \cite{Fer3}
for further discussion). Following \cite{Fer1}, \cite{Fer2}
we introduce the reciprocal invariants:

the symmetric 2-form
$$
\frac{\partial_1\lambda^1\partial_2\lambda^2}{(\lambda^1-\lambda^2)^2}~dR^1dR^2
\eqno(6.9)
$$
and the differential 
$$
d\Omega 
\eqno(6.10)
$$ 
of the 1-form
$$
\Omega=\left(\frac{\partial_1\partial_2\lambda^2}{\partial_2\lambda^2}+
             \frac{\partial_1\lambda^1}{\lambda^1-\lambda^2}\right)dR^1+
       \left(\frac{\partial_1\partial_2\lambda^1}{\partial_1\lambda^1}+
             \frac{\partial_2\lambda^2}{\lambda^2-\lambda^1}\right)dR^2
\eqno(6.11)
$$
($\Omega$ itself is not reciprocally invariant). Note that both objects 
(6.9) and (6.10) do not change under the reparametrization of Riemann invariants
$R^1\to \varphi^1(R^1), ~~ R^2\to \varphi^2(R^2)$.

{\bf Remark}. In order to check the invariance of (6.9) and (6.10) under arbitrary
reciprocal transformations it is sufficient to check their invariance under the 
following elementary ones:
$$
\begin{array}{c}
dX=Bdx+Adt,\\
dT=dt,\\
\end{array} 
$$
which changes only $x$ and preserves $t$ (under this transformation $\lambda^i$
goes to $\Lambda^i=\lambda^iB-A$) and
$$
\begin{array}{c}
dX=dt,\\
dT=dx,\\
\end{array} 
$$
which transforms $\lambda^i$ into $\Lambda^i=\frac{1}{\lambda^i}$.
The invariance of (6.9) and (6.10) under these elementary transformations can
be checked by a direct calculation. Since any reciprocal transformation is a 
composition of elementary ones, we arrive at the required invariance.

It is quite remarkable that invariants (6.9) and (6.10) form a complete 
set in the following sense: if the invariants of one system can be mapped onto the 
invariants of the other one by the  appropriate change of coordinates $R^i$, than both 
these systems are reciprocally related and the corresponding reciprocal transformation
(6.6) can be constructed explicitely
(see \cite{Fer1}, \cite{Fer2} for the discussion).

\section{Hamiltonian systems}

System (6.1) is called Hamiltonian, if it can be represented in the form
$$
u^i_t=\epsilon^i\delta^{ij}\frac{d}{dx}\left(\frac{\delta H}{\delta u^j}\right),
~~ \epsilon^i=\pm 1,
$$
with the Hamiltonian operator  $\epsilon^i\delta^{ij}\frac{d}{dx}$
and the Hamiltonian $H=\int h(u)dx$. In this case the matrix $v^i_j$ is just the 
Hessian of the density $h$ (for definiteness we choose $\epsilon^i=1$), so that system
(6.1) assumes the form
$$
\left(
\begin{array}{c}
u^1\\
u^2
\end{array}
\right)_t=
\left(
\begin{array}{cc}
h_{11} & h_{12}\\
h_{12} & h_{22}
\end{array}
\right)
\left(
\begin{array}{c}
u^1\\
u^2
\end{array}
\right)_x
\eqno(7.1)
$$
(here  $h_{ij}$  means $\frac{\partial^2h}{\partial_{u^i}\partial_{u^j}}$).
For systems (6.2) in Riemann invariants the necessary and sufficient condition for the 
existence of the Hamiltonian representation (7.1) is given by the following 

{\bf Lemma} \cite{Tsarev1}. System (6.2) is Hamiltonian if and only if there exists
flat diagonal metric $ds^2=g_{11}(R){dR^1}^2+g_{22}(R){dR^2}^2$ such that
$$
\begin{array}{c}
\partial_2\ln\sqrt{g_{11}}=\frac{\partial_2\lambda^1}{\lambda^2-\lambda^1},\\
\ \\
\partial_1\ln\sqrt{g_{22}}=\frac{\partial_1\lambda^2}{\lambda^1-\lambda^2}.\\
\end{array}
\eqno(7.2)
$$
Introducing the Lame coefficients $H_1=\sqrt g_{11}, ~~ H_2=\sqrt g_{22}$ and the 
rotation coefficients $\beta_{12}, \beta_{21}$ by the formulae
$$
\partial_1H_2=\beta_{12}H_1, ~~~ \partial_2H_1=\beta_{21}H_2,
\eqno(7.3)
$$
we can rewrite the flatness condition of the metric $ds^2$ in a simple form
$$
\partial_1\beta_{12}+\partial_2\beta_{21}=0.
\eqno(7.4)
$$
The coordinates $u^1, u^2$ in (7.1) are just flat coordinates of the metric $ds^2$, 
where it assumes the standard Euclidean form ${du^1}^2+{du^2}^2$. 

Under  reciprocal transformations (6.6) the metric coefficients $g_{ii}$
transform according to the formulae
$$
G_{ii}=g_{ii}~\frac{(M-\lambda^iN)^2}{(BM-AN)^2}, ~~~~~ i=1, 2
\eqno(7.5)
$$
(see \cite{Serre1}, \cite{Fer4}), so that the transformed metric coefficients
$G_{ii}$ and the transformed characteristic velocities $\Lambda^i$ satisfy the 
same equations (7.2). It is important to emphasize that reciprocal 
transformations do not preserve in general the flatness condition of the metric
$ds^2$ and hence destroy the Hamiltonian structure. However, for any Hamiltonian system
there always exist sufficiently many "canonical" reciprocal transformations
preserving the flatness condition \cite{Fer3}, \cite{Fer4}.

Let us introduce the cubic form
$$
\partial_1\lambda^1g_{11}{dR^1}^3+\partial_2\lambda^2g_{22}{dR^2}^3.
\eqno(7.6)
$$
Using formulae (6.8) and (7.5) one can immediately check, that this qubic form is 
conformally invariant under reciprocal transformations: it aquires the multiple
$\frac{1}{BM-AN}$, so that the zero curves of (7.6) are reciprocally invariant.
Hence  with any Hamiltonian system  we can associate besides the 
invariants (6.9) and (6.10) the reciprocally invariant 3-web of curves
formed by coordinate lines
$R^1=const$, $R^2=const$ and the zero curves of  qubic form
(7.6) which are defined by the equation
$$
(\partial_1\lambda^1g_{11})^{\frac{1}{3}}dR^1+
(\partial_2\lambda^2g_{22})^{\frac{1}{3}}dR^2=0.
$$
A calculation similar to that in sect.4  shows that
 invariant (6.10) is just the curvature form of this 3-web.

{\bf Remark}. It will be interesting to obtain explicit 
formulae for reciprocal invariants (6.9), (6.10) and (7.6) in the flat
coordinates $u^i$ in terms of the Hamiltonian density $h$.

As we already know the objects similar to (6.9), (6.10) and
(7.6) arise in the Lie sphere geometry of surfaces. To clarify this point
we recall the construction of \cite{Fer5}, \cite{Fer3} relating Hamiltonian systems and
surfaces in the Euclidean space.

\section{Hamiltonian systems and surfaces in $E^3$}

Let us consider a 2-component Hamiltonian system (7.1)
$$
\left(
\begin{array}{c}
u^1\\
u^2
\end{array}
\right)_t=
\left(
\begin{array}{cc}
h_{11} & h_{22}\\
h_{12} & h_{22}
\end{array}
\right)
\left(
\begin{array}{c}
u^1\\
u^2
\end{array}
\right)_x
$$
and apply the reciprocal transformation
$$
\begin{array}{c}
dX=Bdx+Adt,\\
dT=dt,\\
\end{array} 
$$
where
$$
B=\frac{{u^1}^2+{u^2}^2+1}{2}, ~~~~~ A=h_1u^1+h_2u^2-h
$$
(this is indeed an integral of system (7.1)). The transformed system assumes the form
$$
\left(
\begin{array}{c}
u^1\\
u^2
\end{array}
\right)_t=
\left(
\begin{array}{cc}
h_{11}B-A & h_{12}B\\
h_{12}B & h_{22}B-A
\end{array}
\right)
\left(
\begin{array}{c}
u^1\\
u^2
\end{array}
\right)_x
\eqno(8.1)
$$
To reveal  geometric meaning of system (8.1) we introduce a surface $M^2$ in the
Euclidean space $E^3(x^1, x^2, x^3)$ with the radius-vector
$$
\vec r=
\left(
\begin{array}{c}
x^1\\
x^2\\
x^3
\end{array}
\right)=
\left(
\begin{array}{c}
h_1-u^1\frac{A}{B}\\
\ \\
h_2-u^2\frac{A}{B}\\
\ \\
-\frac{A}{B}
\end{array}
\right)
\eqno(8.2)
$$
As one can verify by a straightforward calculation, the unit normal of the surface 
$M^2$ is given by  
$$
\vec n=
\left(
\begin{array}{c}
\frac{u^1}{B}\\
\ \\
\frac{u^2}{B}\\
\ \\
\frac{1}{B}-1
\end{array}
\right)
$$
Let us define the matrix $w^i_j$ by the formula
$$
\frac{\partial \vec r}{\partial u^j}=
\sum_{i=1}^2 w^i_j\frac{\partial \vec n}{\partial u^i}.
$$
Geometrically $w^i_j$ is just the inverse of the Weingarten operator (shape operator)
of the surface $M^2$. Using the formulae for $\vec r$ and $\vec n$
we arrive at the following expression for the matrix $w^i_j$:
$$
\left(
\begin{array}{cc}
h_{11}B-A & h_{12}B\\
h_{12}B & h_{22}B-A
\end{array}
\right)
$$
which coincides with that in (8.1). 
Hence the matrix of system (8.1) is just the inverse of the 
Weingarten operator of the associated 
surface $M^2$. The characteristic velocities $w^i$ of system (8.1)
are related to that of system (7.1) by the formula
$$
w^i=\lambda^iB-A,
\eqno(8.3)
$$
and have geometric meaning of radii of principal curvatures of the surface $M^2$.
Moreover, the Riemann invariants of both systems (7.1) and (8.1) coincide and play the
role of parameters of the lines of curvature. 
Equations (8.1) can be equivalently represented in the conservative form
$$
\vec n_t=\vec r_x.
$$

Some further properties of
the correspondence (8.2) (in the general n-component case) were discussed
in \cite{Fer5}, \cite{Fer3}, in particular:

--- commuting Hamiltonian systems correspond via formula (8.2) to surfaces with the 
same spherical image of the lines of curvature;

--- multi-Hamiltonian systems correspond to surfaces, possessing nontrivial deformations
preserving the Weingarten operator;

--- "canonical" reciprocal transformations, preserving the Hamiltonian structure, 
correspond to Lie sphere transformations of the associated surfaces;

--- the flat metric $ds^2$ defining Hamiltonian structure (see lemma in sect.7)
corresponds to the third fundamental form of the associated surface.

Since the correspondence between systems (7.1) and (8.1) is reciprocal, 
invariants (6.9), (6.10) and (7.6) coincide respectively with the symmetric 2-form
$$
\frac{\partial_1w^1\partial_2w^2}{(w^1-w^2)^2}~dR^1dR^2,
$$
the skew-symmetric 2-form $d\Omega$, where
$$
\Omega=\left(\frac{\partial_1\partial_2w^2}{\partial_2w^2}+
             \frac{\partial_1w^1}{w^1-w^2}\right)dR^1+
       \left(\frac{\partial_1\partial_2w^1}{\partial_1w^1}+
             \frac{\partial_2w^2}{w^2-w^1}\right)dR^2
$$
and the conformal class of the qubic form
$$
\partial_1w^1G_{11}{dR^1}^3+\partial_2w^2G_{22}{dR^2}^3
$$
where now $R^1, R^2$ are the parameters of  lines of curvature,
$w^1, w^2$ are the radii of principal
curvatures and $G_{11}, G_{22}$  are the components of the
third fundamental form of the associated surface $M^2$.
Since these objects preserve their form if we rewrite them in terms of principal 
curvatures $k^i$ and the components of the metric $g_{ii}$, they concide with the
Lie sphere invariants of the surface $M^2$. This provides 
remarkable differential-geometric 
interpretation of reciprocal invariants of hydrodynamic type systems.

\section{Reciprocal transformations and reciprocal invariants of 
$n$-component systems}

Let us consider an $n$-component system of hydrodynamic type
$$
u^i_t=v^i_j(u)u^j_x, ~~~~ i, j=1,...,n
\eqno(9.1)
$$
with the characteristic velocities $\lambda^i$ and the corresponding left eigenvectors
$\vec l^i=(l^i_j)$ which satisfy the formulae
$$
\sum_k l^i_kv^k_j=\lambda^il^i_j.
$$
Introducing the 1-forms $\omega^i=l^i_jdu^j$ (note that $\vec l^i$ and $\omega^i$
are defined up to rescailing $\vec l^i\to p^i\vec l^i, ~~ \omega^i \to p^i\omega^i$),
we can rewrite equations (9.1) in the equivalent exterior form
$$
\omega^i \wedge (dx+\lambda^i dt)=0, ~~~~~ i=1,...,n.
\eqno(9.2)
$$
Differentiation of $\omega^i$ and $\lambda^i$ results in the "structure equations" of 
system (9.1):
$$
d\omega^i=c^i_{jk}\omega^j \wedge \omega^k,
\eqno(9.3)
$$
$$
d\lambda^i=\lambda^i_j\omega^j.
\eqno(9.4)
$$
Systems in Riemann invariants are specified by the conditions
$c^i_{jk}=0$ for any triple of indices $i\ne j\ne k$. 
Indeed, in this case the forms $\omega^i$ satisfy the equations
$d\omega^i \wedge \omega^i=0$ for any $i=1,...,n$ and hence can be normalized 
so as to become just $\omega^i=dR^i$. In the coordinates $R^i$ equations (9.2)
assume the familiar Riemann-invariant form
$$
R^i_t=\lambda^i R^i_x, ~~~~~ i=1,...,n.
\eqno(9.5)
$$
The exterior representation (9.2) is a natural analog of  representation (9.5) 
which is applicable in the nondiagonalizable case as well. We emphasize that for 
$n\geq 3$ Riemann invariants do not exist in general. 

Applying to (9.1) the
reciprocal transformation
$$
\begin{array}{c}
dX=Bdx+Adt,\\
dT=Ndx+Mdt,\\
\end{array} 
$$
we arrive at the transformed equations
$$
u_T^i=V^i_j(u)u^j_X
$$
with the new matrix $V$ given by
$$
V=(Bv-AE)(ME-Nv)^{-1}, ~~~~ E=id
$$
or, in the exterior form,
$$
\omega^i \wedge (dX+\Lambda^i dT)=0
$$
where
$$
\Lambda^i=\frac{\lambda^iB-A}{M-\lambda^iN}.
$$
Hence the forms $\omega^i$ as well as the 
structure equations (9.3) do not change, while $\lambda^i$ transform as in the 
2-component case --- see formula (6.8).

{\bf Remark}. In the $n$-component case equations (6.4) for the densities and fluxes of 
conservation laws $hdx+gdt$ assume the form
$$
g_i=\lambda^ih_i, ~~~~ i=1,...,n
$$
where $g_i$ and $h_i$ are defined by the expansions
$$
dg=g_i\omega^i, ~~~ dh=h_i\omega^i.
$$

In \cite{Fer1}, \cite{Fer2} we introduced the following reciprocally invariant objects:

1. The symmetric 2-form
$$
\sum_{i\ne j}\frac{\lambda^i_i\lambda^j_j}{(\lambda^i-\lambda^j)^2}~\omega^i\omega^j
\eqno(9.6)
$$

2. The skew-symmetric 2-form 
$$
d\Omega
\eqno(9.7)
$$
where
$$
\Omega=\left(\sum_{k\ne 1}\frac{\lambda^k_1-\frac{\lambda^1_1}{n-1}}{\lambda^1-\lambda^k}
         \right)\omega^1+\ldots +
       \left(\sum_{k\ne n}\frac{\lambda^k_n-\frac{\lambda^n_n}{n-1}}{\lambda^n-\lambda^k}
         \right)\omega^n
\eqno(9.8)
$$
($\Omega$ itself is not reciprocally invariant). Note that both objects (9.6) and (9.7)
do not change if we reparametrize the 1-forms in the structure equations:
$\omega^i\to p^i\omega^i$. Objects (9.6) and (9.7) are the natural analogs of the
corresponding invariants (6.9) and (6.10) in the 2-component case. However, for 
$n\geq 3$ the form $\Omega$ depends only on the first derivatives of the characteristic 
velocities $\lambda^i$ rather that on the second derivatives as in the 2-component case.

In principle for $n\geq 3$ there exist additional reciprocal invariants, namely
the 1-forms
$$
\frac{\lambda^i_i(\lambda^j-\lambda^l)}{(\lambda^i-\lambda^j)(\lambda^i-\lambda^l)}
\omega^i ~~~~
(i\ne j\ne l)
$$
as well as the cross-ratios
$$
\frac{(\lambda^i-\lambda^j)(\lambda^k-\lambda^l)}
     {(\lambda^k-\lambda^j)(\lambda^i-\lambda^l)}
$$
of any four characteristic velocities. 

However, as follows from 
\cite{Fer1}, \cite{Fer2}, the structure equations (9.3) and the invariants (9.6), (9.7)
 in fact define generic system of hydrodynamic type uniquely up to reciprocal 
 equivalence (under "generic" it is sufficient to understand genuinely nonlinear
  system, that is, a system with $\lambda^i_i\ne 0$ for any $i$).

\section{Appendix. Derivation of the mVN equation for diagonally-cyclidic surfaces}

Our aim is to show, that system (4.4) on $\rho, k^1, k^2$ can be rewritten in
terms of the Lie sphere density $U$ defined as
$$
U^2=\frac{\partial_1k^1\partial_2k^2}{(k^1-k^2)^2}.
$$
In view of the correspondence between surfaces in $E^3$ and Hamiltonian
systems discussed in sect. 8 the classification of diagonally-cyclidic surfaces can
be equivalently reformulated as the classification of those Hamiltonian
systems
$$
\begin{array}{c}
R_t^1=\lambda^1R^1_x \\
R_t^2=\lambda^2R^2_x \\
\end{array} 
$$
for which the flat metric $ds^2=g_{11}{dR^1}^2+g_{22}{dR^2}^2$ defining the
Hamiltonian structure can be represented in the form
$$
g_{11}=\frac{e^{2 \rho}}{\partial_1\lambda^1}, ~~~~
g_{22}=\frac{e^{2 \rho}}{\partial_2\lambda^2}.
$$
Indeed, this is an immediate consequence of isotermicity of the cubic form.
Equations (7.2) and the flatness condition for the metric $ds^2$ imply the
following system for $\rho, \lambda^1, \lambda^2$:
$$
\begin{array}{c}
\partial_1  \rho = B+\frac{1}{2}\frac{\partial_1\partial_2\lambda^2}
{\partial_2\lambda^2}, \\
\ \\
\partial_2  \rho = A+\frac{1}{2}\frac{\partial_1\partial_2\lambda^1}
{\partial_1\lambda^1}, \\
\end{array}
\eqno(10.1)
$$
$$
\partial_1\lambda^1\left(\partial_1B+\frac{B}{2}\partial_1\ln
{\frac{\partial_1\lambda^1}{\partial_2\lambda^2}}\right)+
\partial_2\lambda^2\left(\partial_2A+\frac{A}{2}\partial_2\ln
{\frac{\partial_2\lambda^2}{\partial_1\lambda^1}}\right)=0
\eqno(10.2)
$$
where
$A=\frac{\partial_2\lambda^1}{\lambda^2-\lambda^1}, ~~
B=\frac{\partial_1\lambda^2}{\lambda^1-\lambda^2}$ (compare with (4.4)).
Since
$$
U^2=\frac{\partial_1k^1\partial_2k^2}{(k^1-k^2)^2}=
\frac{\partial_1\lambda^1\partial_2\lambda^2}{(\lambda^1-\lambda^2)^2}
$$
the desired equation for $U$ will be obtained after we rewrite (10.1) --
(10.2) in terms of $U$. Crossdifferentiating (10.1) and introducing
$k=\frac{1}{2}\ln{\frac{\partial_2\lambda^2}{\partial_1\lambda^1}}$
we obtain
$$
\partial_1\partial_2k=\partial_1A-\partial_2B,
$$
while (10.2) assumes the form
$$
\partial_2(e^kA)+\partial_1(e^{-k}B)=0.
$$
Moreover, the equations for the characteristic
velocities
$$
\partial_1\lambda_1=Ue^{-k}(\lambda^1-\lambda^2), ~~~~
\partial_2\lambda_1=A(\lambda^2-\lambda^1)
$$
$$
\partial_1\lambda_2=B(\lambda^1-\lambda^2), ~~~~
\partial_2\lambda_2=Ue^{k}(\lambda^1-\lambda^2)
$$
give as the compatibility conditions two additional relations
$$
\begin{array}{c}
\partial_1A=AB+U^2-\partial_2(Ue^{-k}), \\
\ \\
\partial_2B=AB+U^2+\partial_1(Ue^{k}).
\end{array}
$$
Thus the problem is reduced to  that of excluding variables $A, B, k$ and
deriving the equation for $U$ from the following much simpler looking system
$$
\partial_1\partial_2k=\partial_1A-\partial_2B, \eqno(10.3)
$$
$$
\partial_2(e^kA)+\partial_1(e^{-k}B)=0, \eqno(10.4)
$$
$$
\partial_1A=AB+U^2-\partial_2(Ue^{-k}), \eqno(10.5)
$$
$$
\partial_2B=AB+U^2+\partial_1(Ue^{k}).
\eqno(10.6)
$$
To proceed further we introduce new variables $m, n$ by the formulae
$$
-\partial_2(Ue^{-k})=Um, ~~~~~ \partial_1(Ue^{k})=Un,
$$
so that $\partial_1k, \partial_2k$ can be expressed as follows:
$$
\partial_1k=-\frac{\partial_1U}{U}+e^{-k}n, ~~~~~
\partial_2k=\frac{\partial_2U}{U}+e^{k}m.
\eqno(10.7)
$$
The compatibility conditions of (10.7) with (10.3) (which now assumes the form
$\partial_1\partial_2k=U(m-n)$) give  the following equations for $m, n$:
$$
\begin{array}{c}
\partial_2n=\frac{\partial_2U}{U}n+e^{k}(\partial_1\partial_2\ln U+U(m-n)+mn) \\
\ \\
\partial_1m=\frac{\partial_1U}{U}m-e^{-k}(\partial_1\partial_2\ln U-U(m-n)+mn). \\
\end{array}
\eqno(10.8)
$$
Rewriting (10.4) in the form
$$
e^k\partial_2A+e^kA(\frac{\partial_2U}{U}+e^km)+
e^{-k}\partial_1B-e^{-k}B(-\frac{\partial_1U}{U}+e^{-k}n)=0
$$
and introducing the new variable $F$ by the formulae
$$
e^k\partial_2A+e^kA(\frac{\partial_2U}{U}+e^km)=
F+\frac{1}{2}e^{k}A^2-\frac{1}{2}e^{-k}B^2,
$$
$$
e^{-k}\partial_1B-e^{-k}B(-\frac{\partial_1U}{U}+e^{-k}n)=
-F-\frac{1}{2}e^{k}A^2+\frac{1}{2}e^{-k}B^2,
$$
we can express the derivatives of $A$ and $B$ as follows
$$
\begin{array}{c}
\partial_1A=AB+U^2+Um, \\
\ \\ 
\partial_2B=AB+U^2+Un, \\
\ \\
\partial_2A=-A(\frac{\partial_2U}{U}+e^km)+
\frac{1}{2}A^2-\frac{1}{2}e^{-2k}B^2+e^{-k}F, \\
\ \\
\partial_1B=-B(\frac{\partial_1U}{U}-e^{-k}n)
-\frac{1}{2}e^{2k}A^2+\frac{1}{2}B^2-e^{k}F. \\
\end{array}
\eqno(10.9)
$$
Let us introduce also the new functions $G$ and $H$ by the formulae
$$
\begin{array}{c}
\partial_1n=e^{2k}F+\frac{1}{2}e^{-k}(3n^2+4nU+U^2)-2n\frac{\partial_1U}{U}-
3\partial_1U+\frac{1}{2}e^{k}H, \\
\ \\
\partial_2m=e^{-2k}F-\frac{1}{2}e^{k}(3m^2+4mU+U^2)-2m\frac{\partial_2U}{U}-
3\partial_2U+\frac{1}{2}e^{-k}G. \\
\end{array}
\eqno(10.10)
$$
The compatibility conditions of eqns. (10.9) give the following expressions for 
$\partial_1F, \partial_2F$:
$$
\begin{array}{c}
\partial_1F=- F(\frac{\partial_1U}{U}-e^{-k}n)+e^{2k}Um(U+m) \\
\ \\
+\frac{U}{2}(G+2e^{-k}F-e^{2k}(3m^2+4mU+U^2)), \\
\ \\
\partial_2F=- F(\frac{\partial_2U}{U}+e^{k}m)+e^{-2k}Un(U+n) \\
\ \\
-\frac{U}{2}(H+2e^{k}F+e^{-2k}(3n^2+4nU+U^2)). \\
\end{array}
\eqno(10.11)
$$
The compatibility conditions of (10.8) and (10.10) give the expressions for 
$\partial_1G, \partial_2H$ which can be represented as follows:
$$
\begin{array}{c}
\partial_1(G+2\frac{\partial_2^2U}{U}-(\frac{\partial_2U}{U})^2)=-3\partial_2(U^2), \\
\ \\
\partial_2(H-2\frac{\partial_1^2U}{U}+(\frac{\partial_1U}{U})^2)=3\partial_1(U^2). \\
\end{array}
\eqno(10.12)
$$
Finally, the compatibility conditions of (10.11) imply
$$
\partial_2(GU^2)+\partial_1(HU^2)=0,
\eqno(10.13)
$$
so that equations (10.3)-(10.6) are reduced to (10.12)-(10.13).
The desired stationary mVN equation results now after introducing $V, W$ by the
formulae
$$
G+2\frac{\partial_2^2U}{U}-\left(\frac{\partial_2U}{U}\right)^2=-3W,
$$
$$
H-2\frac{\partial_1^2U}{U}+\left(\frac{\partial_1U}{U}\right)^2=3V.
$$
All these calculations were performed with Mathematica
and can be easily verified.

\section{Concluding remarks}

Here we just list some of the unsolved problems.

1. In our discussion of Lie-sphere invariants of surfaces
(reciprocal invariants of hydrodynamic-type systems)
the choice of coordinates $R^i$ plays a crucial
role. In case of surfaces these are coordinates of the lines of curvature
(Riemann invariants in case of hydrodynamic type systems).
This choice
is not accidental, since the lines of curvature are preserved by the Lie sphere
group while Riemann invariants are preserved under
reciprocal transformations.
In fact only in these special coordinates do our invariants assume
particularly symmetric and simple form. However, from the point of view of
applications it is desirable to have a kind of invariant tensor formula, which
will allow computation of these objects in an arbitrary coordinate system,
for instance, in conformal parametrization in case of surfaces or
in flat coordinates in case of Hamiltonian systems.
   
2. In \cite{Kon1}, \cite{Kon2} Konopelchenko introduced dynamics of surfaces,
governed by the mVN equation. In this approach the
integrals of mVN correspond to certain functionals on surfaces, which 
were conjectured in \cite{Taiman3}, \cite{Taiman1}, 
\cite{Taiman2} to be conformally invariant. In particular,
the simplest quadratic integral of mVN  corresponds to the
Willmore functional. This conjecture was proved recently in
\cite{Grin}. Since functional (1.3) is conformally invariant
(indeed, it is invariant under the full group of Lie sphere 
transformations which contains conformal group), 
it would be interesting to understand its relationship to the mVN hierarchy.

3. It seems to be an interesting problem to describe the class of surfaces,
  for which evolutions (3.5), (3.6) preserve the lines of curvature.

\section{Acknowledgements}

I would like to thank  A.I.~Bobenko, B.G.~Konopelchenko, 
U.~Pinkall and colleagues in 
Technische Universit\"at, Berlin for their interest and useful remarks.
I am especially grateful to K.R. Khusnutdinova for performing a number of important
calculations with Mathematica without which  the paper
would not appear in its present form.

This research was supported by  the Aleksander von Humboldt Foundation.


\begin{thebibliography}{99}
\addcontentsline{toc}{section}{References}


\bibitem{Lie} Lie~S., \"Uber Komplexe, inbesondere Linien- und Kugelkomplexe,
mit Anwendung auf der Theorie der partieller Differentialgleichungen,
Math. Annalen, 1872, V.5, 145-208, 209-256.

\bibitem{Bla3} Blaschke~W., Vorlesungen \"uber Differentialgeometrie, V.3,
Springer-Verlag, Berlin, 1929.


\bibitem{Pinkall} Pinkall~U., Dupinische Hyperfl\"achen in $E^4$, 
Manuscripta Math., 1985, V.51, 89-119.

\bibitem{Pinkall2} Pinkall~U., Dupin hypersurfaces, Math. Annalen,
1985, V.270, 427-440.


\bibitem{Cecil2} Cecil~T., Lie sphere geometry, Springer-Verlag, 1992.

\bibitem{Sym} Cieslinski~J., Goldstein~P. and Sym~A.,
 Isothermic surfaces in $E^3$ as soliton surfaces, Phys. lett. A,
1995, V.205, 37-43.

\bibitem{Bobenko} Bobenko~A. and Pinkall~U., Discrete isothermic surfaces,
J. Reine Angew. Math., 1996, V.475, 187-208.

\bibitem{Calapso} Calapso~P., Sulla superficie a linee di curvatura isoterme,
Rend. Circ. Mat. Palermo, 1903, V.17, 275-286.

\bibitem{Kon1} Konopelchenko~B.G., Multidimensional integrable systems and dynamics of
surfaces in space, Preprint Institute of Math., Acad. Sinica, Taipei, Taiwan, R.O.C.,
1993. 

\bibitem{Kon2} Konopelchenko~B.G., Induced surfaces and their integrable dynamics,
Studies in Appl. Math., 1996, 9-51.

\bibitem{Kon3} Konopelchenko~B.G., Nets in $R^3$, their integrable evolutions and 
the DS hierarchy, Phys. letters A, 1993, V.183, 153-159.

\bibitem{Taiman3} Taimanov~I.A., Modified Novikov-Veselov equation and differential 
geometry of surfaces, in Solitons, Geometry and Topology 
(eds. V.M.Buchstaber and S.P.Novikov) Transl. AMS, ser.2, 1997, V.179, 133-155.

\bibitem{Taiman1} Taimanov~I.A., Surfaces of revolution in terms of solitons, 
Preprint  (dg-ga/9610013) to appear in Ann. of Global Anal. and Geometry,
1997, V.15, N5.

\bibitem{Taiman2} Taimanov~I.A., Global Weierstrass representation and it's spectrum,
to appear in Russian Math. Surveys, 1997.




\bibitem{Grin} Grinevich~P.G. and Schmidt~M.U., Conformally invariant
functionals of immersions of tori into $R^3$, Preprint SFB 288 N 252, Berlin, 1997.









\bibitem{Rogers1} Rogers~C. and Shadwick~W.F., B\"acklund transformations 
and their applications, Academic Press, New York, 1982.

\bibitem{Rogers2} Rogers~C., Reciprocal transformations and their applications,
Nonlinear Evolutions, Proc. of the 5th Workshop on Nonlinear Evolution Equations 
and Dynamical systems, France, 1987, 109-123.


\bibitem{Fer1} Ferapontov~E.V., Reciprocal transformations and their invariants,
Diff. Uravneniya, 1989, V.25, N.7, 1256-1265 (English translation in
Differential Equations, 1989, V.25, N.7, 898-905).

\bibitem{Fer2} Ferapontov~E.V., Reciprocal autotransformations and hydrodynamic 
symmetries, Diff. Uravneniya, 1991, V.27, N.7, 1250-1263 (English translation in
Differential Equations, 1989, V.27, N.7, 885-895). 

\bibitem{Fer5} Ferapontov~E.V., Hamiltonian systems of hydrodynamic type and their 
realization on hypersurfaces of a pseudoeuclidean space, Soviet J. Math., 1991,
V.55, 1970-1995.

\bibitem{Fer3} Ferapontov~E.V., Dupin hypersurfaces and integrable Hamiltonian 
systems of hydrodynamic type which do not possess Riemann invariants, Diff. Geometry
and its Appl., 1995, V.5, 121-152.




\bibitem{Serre1} Serre~D., Oscillations non lineaires des systemes hyperboliques:
methodes et resultats qualitatifs, Ann. Inst. Henri Poincare,
1991, V.8, N.3-4, 351-417.



\bibitem{Fer4} Ferapontov~E.V., Nonlocal Hamiltonian operators of hydrodynamic type:
Differential geometry and Applications, Amer. Math. Soc. Transl., 
1995, (2) v.170, 33-58.


\bibitem{Fer6} Ferapontov~E.V., On integrability of $3\times 3$ semihamiltonian 
hydrodynamic  type systems  which do not possess Riemann invariants, Physica D,
1993, V.63, 50-70.


\bibitem{Tsarev1} Tsarev~S.P., The geometry of Hamiltonian systems of hydrodynamic type.
The generalized hodograph transform, Math. USSR Izv., 1991, V.37, 397-419.


\bibitem{Bla1} Blaschke~W. and Bol~G., Geometrie der Gewebe, Springer-Verlag, 
Berlin, 1938.

\bibitem{Bla2} Blaschke~W., Einf\"uhrung in die Geometrie der Waben, 
Birkh\"auser-Verlag, Basel, Switzerland, 1955.



\bibitem{Wang} Changping Wang, Surfaces in M\"obius geometry, Nagoya Math. J.,
1992, V.125, 53-72.


\bibitem{AkGold1} Akivis~M.A. and Goldberg~V.V.,
Conformal differential geometry and it's generalizations,
New York, Wiley, 1996.

\bibitem{AkGold2} Akivis~M.A. and Goldberg~V.V.,
Projective differential geometry of submanifolds, Math. Library, V.49,
North-Holland, 1993.


\bibitem{Mikh} Mikhailov~A.V. and Yamilov~R.I. On integrable two-dimensional
generalizations of Nonlinear Schr\"odinger type Equations, 
Phys. letters A. 1997, V.230, 295-300.

\bibitem{Cecil1} Cecil~T., On the Lie curvature of Dupin hypersurfaces,
Kodai Math. J., 1990,  V.13, 143-153.



\bibitem{Miyaoka} Miyaoka~R., Dupin hypersurfaces and a Lie invariant,
Kodai Math. J., 1989,  V.12, 228-256.



\bibitem{Bogdanov} Bogdanov~L.V., Veselov-Novikov equation as a natural 
two-dimensional generalization of the Korteweg-de Vries equation, Theor. 
and Math. Phys., 1987, V.70, 309-314.

\bibitem{Bobenko2} Bobenko~A. and Schief~W., Discrete affine spheres,
Preprint SFB 288 No.263, Berlin, 1997. 


\bibitem{Finikov} Finikov S.P., Theory of congruences, Moscow-Leningrad, 1950.

\bibitem{Kon4} Konopelchenko~B.G. and Rogers~C., On generalized Loewner system:
novel integrable equations in 2+1-dimensions, J. Math. Phys., 1993, V.34, N1,
214-242.

\end{thebibliography}
\end{document}